\documentstyle[12pt,amssymb,amscd]{amsart}  % amslatex

\newtheorem{thm}{Theorem}[section]
\newtheorem{cor}[thm]{Corollary}
\newtheorem{prop}[thm]{Proposition}
\newtheorem{lemma}[thm]{Lemma}

\theoremstyle{remark}
\newtheorem{remark}[thm]{Remark}
\newtheorem{example}[thm]{Example}

\theoremstyle{definition}
\newtheorem{defn}[thm]{Definition}

\numberwithin{equation}{section}

\newcommand{\bbC}{{\Bbb C}}

\newcommand{\bbR}{{\Bbb R}}

\newcommand{\bbZ}{{\Bbb Z}}

\newcommand{\cF}{{\cal F}}

\newcommand{\cO}{{\cal O}}

\newcommand{\cL}{{\cal L}}
\newcommand{\cM}{{\cal M}}
\newcommand{\cN}{{\cal N}}

\newcommand{\cV}{{\cal V}}
\newcommand{\cA}{{\cal A}}
\newcommand{\cB}{{\cal B}}

\newcommand{\cC}{{\cal C}}
\newcommand{\cE}{{\cal E}}
\newcommand{\cW}{{\cal W}}
\newcommand{\cU}{{\cal U}}

\newcommand{\cK}{{\cal K}}
\newcommand{\cH}{{\cal H}}

\newcommand{\Hom}{\operatorname{Hom}}
\newcommand{\Ext}{\operatorname{Ext}}

\newcommand{\de}{\operatorname{def}}
\newcommand{\Tot}{\operatorname{Tot}}
\newcommand{\End}{\operatorname{End}}

\newcommand{\printname}[1]
  {\smash{\makebox[0pt]{\hspace{-2.0in}\raisebox{8pt}{\tiny #1}}}}

\thanks{This research was partially supported by the CRDF grant RM1-2089 
and by the NSA grant MDA904-01-1-0020.}

\begin{document}

\title{Deformations of quasicoherent sheaves of algebras}

\author{Valery A.~Lunts}
\address{Department of Mathematics, Indiana University,
Bloomington, IN 47405, USA}
\email{vlunts@@indiana.edu}

\begin{abstract} Gerstenhaber and Schack ([GS]) developed a deformation 
theory of presheaves of algebras on small categories. We translate their 
cohomological description to sheaf cohomology. More precisely, we describe 
the deformation space of (admissible) quasicoherent sheaves of algebras 
on a quasiprojective scheme $X$ in terms of sheaf cohomology on $X$ and 
$X\times X$. These results are applied to the study of deformations of the 
sheaf $D_X$ of differential operators on $X$. In particular, in case $X$ 
is a flag variety we show that any deformation of $D_X$, which is induced 
by a deformation of $\cO_X$, must be trivial. This result is used in [LR3], 
where we study the localization construction for quantum groups.
\end{abstract}

\maketitle

\section{Introduction}

Let $X$ be a topological space, $k$ be a field, and $\cA _X$ be a sheaf 
of $k$-algebras on $X$. We would like to study infinitesimal deformations of 
$\cA _X$. Such deformatioms form a $k$-vector space which we denote by 
$\de (\cA _X)$. In case $X=pt$ it is well known that the infinitesimal 
deformations of (the $k$-algebra) $A=\cA _X$ are controlled by the Hochschild 
cohomology of $A$. More precisely, $\de (A)=HH^2(A)=\Ext ^2_{A\otimes A^o}
(A,A)$. However, for a general $X$ and $\cA _X$ the situation is more subtle. 
More generally, 
given an $\cA _X$-bimodule $\cM _X$ we may ask for cohomological 
interpretation of $exal(\cA _X,\cM _X)$ -- the space of algebra 
extensions of 
$\cA _X$ by $\cM _X$ ($exal(\cA _X,\cA _X)=
\de (\cA _X)$). 

Gerstenhaber and Schack ([GS]) developed a deformation theory of 
{\it presheaves} 
of algebras. Given a small category $\cU $ and a presheaf of algebras 
$\cA _{\cU}$ 
on $\cU$ (i.e. a contravariant functor from $\cU$ to the category of 
$k$-algebras) they consider the space $\de (\cA _{\cU})$ of infinitesimal 
deformations of $\cA _{\cU}$ and give it a cohomological interpretation. 
Namely, given an $\cA _{\cU}$-bimodule $\cM _{\cU }$ they define a natural 
exact sequence of complexes of $k$-vector spaces
$$0\to T_a^\bullet(\cM _{\cU})\to T^\bullet (\cM _{\cU})
\to \bar{T}^\bullet (\cM _{\cU})\to 0.$$
The middle term is the total complex of the {\it simplicial bar 
resolution} of $\cM _{\cU}$ and 
$$H^i(T ^\bullet(\cM _{\cU}))=
\Ext^i_{\cA _{\cU}\otimes \cA _{\cU}^o}(\cA _{\cU},\cM _{\cU})$$
-- the Hochschild cohomology of $\cA _{\cU}$ with coefficients in $\cM _{\cU}$. 
The cohomology $H^i(\bar{T}^\bullet (\cM _{\cU}))$ is the cohomology 
$H^i(\cU,\cM_{\cU})$ of 
the nerve of $\cU$ (or the classifying space of $\cU$) 
with coefficients in $\cM _{\cU}$. Finally, 
$$H^2(T_a ^\bullet (\cM _{\cU}))=exal(\cA _{\cU},\cM _{\cU});$$ 
in particular, $H^2(T_a ^\bullet (\cA _{\cU }))=
\de(\cA _{\cU})$. As a consequence they obtain a long exact sequence of $k$-
spaces
$$\begin{array}{cccccc}
  ... & \to & \Ext ^1_{\cA _{\cU}\otimes \cA _{\cU }^o}(\cA _{\cU},
\cM _{\cU}) & \to & H^1(\cU ,\cM _{\cU}) & \to \\
exal(\cA _{\cU},\cM _{\cU}) & \to & \Ext ^2_{\cA _{\cU}\otimes \cA _{\cU }^o}
(\cA _{\cU}, \cM _{\cU}) & \to & H^2(\cU ,\cM _{\cU}) & \to \\
 ... & & & & & \\
\end{array}$$

Returning to our problem of trying to interpret cohomologically the space 
$exal(\cA _X,\cM _X)$ we may proceed as follows. Let $\cU$ be the category 
of (all or some) open subsets of $X$. From the sheaf of algebras $\cA _X$ and 
its bimodule $\cM _X$ we obtain the corresponding presheaves $\cA _{\cU}$ and 
$\cM _{\cU}$. At this point there are two natural questions.

\medskip
\noindent {\bf Q1}. Is $exal(\cA _X,\cM _X)$ equal to $exal(\cA _{\cU},\cM _{\cU})$? 

\noindent {\bf Q2}. Can we interpret the spaces  
$\Ext^i_{\cA _{\cU}\otimes \cA _{\cU}^o}(\cA _{\cU},\cM _{\cU})$ 
 and $H^i(\cU, \cM _{\cU})$ as 
sheaf cohomologies on $X$ or $X\times X$?  

\medskip

The answers to these questions in general are probably negative. 

\medskip

In this paper we obtain positive answers to the above questions in case 
$X$ is a quasiprojective scheme over $k$ and $\cA_X$ and $\cM_X$ are 
quasicoherent sheaves on $X$, which satisfy some additional conditions (the 
pair $(\cA_X,\cM_X)$ must be admissible in the sense of Definition  4.7 
below). In this case there is a natural quasicoherent sheaf of algebras 
$\cA^e_Y$ on the product scheme $Y=X\times X$ (this is the analogue of the 
ring $A\otimes A^o$ for a single algebra $A$). Moreover, the $\cA_X$-bimodule 
$\cM_X$ gives rise to a $\cA^e_Y$-module $\tilde{\cM}_Y$; 
in particular, the $\cA_X$-bimodule $\cA_X$ defines an 
$\cA^e_Y$-module $\tilde{\cA}_Y$. If $\cU$ is the category 
of all {\it affine} open subsets of $X$, then we prove that 
$$exal(\cA_X,\cM_X)=exal(\cA_{\cU},\cM_{\cU}),$$
and
$$\Ext^i_{\cA _{\cU}\otimes \cA _{\cU}^o}(\cA _{\cU},\cM _{\cU})=
\Ext^i_{\cA^e_Y}(\tilde{\cA}_Y, \tilde{\cM}_Y),$$
$$H^i(\cU, \cM _{\cU})=H^i(X,\cM_X).$$

In particular, we obtain the long exact 
sequence 
$$\begin{array}{cccccc}
  ... & \to & \Ext ^1_{\cA^e_Y}(\tilde{\cA} _Y,
\tilde{\cM} _Y) & \to & H^1(X ,\cM _X) &\to  \\
 exal(\cA _X,\cM _X) & \to & \Ext ^2_{\cA^e_Y}
(\tilde{\cA} _Y, \tilde{\cM} _Y) & \to & H^2(X ,\cM _{X}) &\to \\
 ... & & & & & \\
\end{array}$$
which allows us to analyze the space $exal(\cA _X,\cM _X)$. One of the 
implications is that $exal(\cA _X,\cM _X)$ behaves well with respect to base 
field extensions. It is easy to describe the morphisms 
$$H^1(X,\cM _X)\to exal(\cA _X,\cM _X)\to \Ext ^2_{\cA^e_Y}
(\tilde{\cA} _Y,\tilde{\cM} _Y)$$
explicitly. Note that if $X$ is affine then 
$H^i(X,\cM_X)=0$ for $i>0$ and hence $exal(\cA_X,\cM_X)=\Ext ^2
_{\cA^e_Y}(\tilde{\cA}_Y,\tilde{\cM}_Y)$. Moreover, in this case 
$$\Ext ^\bullet
_{\cA^e_Y}(\tilde{\cA}_Y,\tilde{\cM}_Y)=\Ext ^\bullet
_{\cA_X(X)\otimes \cA_X^o(X)}(\cA_X(X),\cM_X(X))$$
and thus 
$$exal(\cA_X,\cM_X)=exal(\cA_X(X),\cM_X(X)).$$

In the special case when $\cA _X=\cO _X$ and 
$\cM _X$ is a symmetric $\cO _X$-bimodule  
the isomorphism 
$$\Ext ^i_{\cO_Y}(\cO _X,\cM _X)=
\Ext ^i_{\cA _{\cU}\otimes \cA _{\cU}^o}(\cA _{\cU},\cM _{\cU})$$
was proved by R.~Swan in [S].

We apply the above results to analyze $\de (\cA _X)$ in case $X$ is 
 a smooth quasiprojective variety over $\bbC$ and $\cA _X=D_X$ -- the sheaf 
of differential operators on $X$. 
In this case 
$$\Ext ^i_{\cA^e_Y}(\tilde{\cA}_Y,\tilde{\cA}_Y)=H^i(X^{an},\bbC ).$$ 
 
If in addition $X$ is $D$-affine (for example $X$ is affine) 
then $H^i(X,D_X)=0$ for $i>0$ and hence 
$$\de(D_X)=H^2(X^{an},\bbC ).$$

In the last section we study {\it induced} deformations of $D_X$, 
i.e. those which come 
from deformations of the structure sheaf $\cO _X$. In particular if $X$ is 
a flag variety we show that every induced deformation of $D_X$ is trivial. 
This result is used in the work [LR3], where we study 
quantum differential operators on quantum flag varieties.

It is my pleasure to thank Paul Bressler for his references to 
the literature on the deformation theory and Michael Larsen for helpful 
discussions of the subject.   

\section{Preliminaries on extension of algebras and Hochschild cohomology}

\subsection{Extensions of algebras}
Fix a field $k$. An algebra means an associative unital $k$-algebra. 
Fix an algebra $A$;   
$A^o$ is the opposite algebra and $A^e:=A\otimes _kA^o$. An $A$-module 
means a left $A$-module; an $A$-bimodule means an $A^e$-module.

Fix an algebra $A$ and an $A$-bimodule $M$. Consider an exact sequence 
of $k$-modules 
$$0\to M\to B\stackrel{\epsilon}{\to}A\to 0$$
with the following properties

\begin{itemize}
\item $B$ is an algebra and $\epsilon$ is a homomorphism of algebras. 
(Hence $M$ is a 2-sided ideal in $B$.)

\item The $B$-bimodule structure on $M$ factors through the homomorphism 
$\epsilon$ and the resulting $A$-bimodule structure on $M$ coincides with 
the given one. (In particular, the square of the ideal $M$ is zero.)
\end{itemize}

\begin{defn} An exact sequence as above is called an {\it algebra extension} 
of $A$ by $M$. An isomorphism between extensions 
$$0\to M\to B\to A\to 0$$
and 
$$0\to M\to B^\prime \to A\to 0$$
is an isomorphism of algebras $\alpha :B\to B^\prime$ which makes the 
 following diagram commutative
$$
\begin{array}{ccrcrcrcc}
0 & \to & M & \to & B & \to & A & \to & 0 \\
  &     & id \downarrow & & \alpha \downarrow & & id \downarrow & & \\
0 & \to & M & \to & B^\prime & \to & A & \to & 0
\end{array}
$$
An extension is {\it split} if there exists an algebra homomorphism 
$s:A\to B$ such that $\epsilon \cdot s=id$. Then $B=A\oplus M$ with 
the multiplication $(a,m)(a^\prime ,m^\prime )=(aa^\prime ,am^\prime +
ma^\prime )$. The collection of isomorphism classes of algebra extensions 
of $A$ by $M$ is naturally a $k$-vector space which is denoted $exal(A,M)$. 
The zero element is the class of the split extension.
\end{defn}

Given a map of $A$-bimodules $M\to M^\prime$ the usual pushout construction 
for extensions defines a map 
$$exal(A,M)\to exal(A,M^\prime).$$
Given a homomorphism of algebras $A^\prime \to A$ the pullback construction 
for extensions defines a map 
$$exal(A,M)\to exal (A^\prime ,M).$$
Thus $exal(\cdot,\cdot)$ is a bifunctor covariant in the second variable and 
contravariant in the first one. 

In case $M=A$ the space $exal(A,A)$ can be considered as deformations of the 
first order of the algebra $A$. Let us describe this space in a different 
way. Put $k_1:=k[t]/(t^2)$. Consider $k_1$-algebras $B$ with a given 
isomorphism $\theta :grB\to A\otimes _kk_1$. (The algebra $B$ has the 
filtration $\{0\}\subset tB \subset B$ and $grB$ denotes the associated 
graded.) The isomorphism classes of such pairs $(B,\theta )$ form a pointed 
set  
which we denote by $\de (A)$. The distinguished element in $\de (A)$ 
is represented by the algebra 
$B=A\otimes _kk_1$. 

We claim that $exal(A,A)=\de (A)$ (hence $\de (A)$ is a $k$-vector space). 
Indeed, given $(B,\theta )$ as above 
we obtain an exact sequence 
$$0\to tB=A\to B\to A\to 0,$$ 
which gives a well defined map from $\de (A)\to exal(A,A)$.
Conversely, given an algebra extension 
$$0\to M=A\to B\to A\to 0$$
define the multiplication $t:B\to B$ by $t\cdot 1_B=1_A\in M$. This makes 
$B$ a $k_1$-algebra and defines the inverse map $exal(A,A)\to \de (A)$. 

The above description of $exal(A,A)$ allows us to define the set  
$\de ^n(A)$  of $n$-th order deformations of $A$ as the collection of 
isomorphism classes of $k_n:=k[t]/(t^{n+1})$-algebras $B$ with an isomorphism 
of  $k_n$-algebras $grB\to A\otimes _kk_n$. 
Thus $\de ^1(A)=\de (A)=exal(A,A)$. The algebra $B=A\otimes _kk_n$ represents 
the  {\it trivial} deformation. Note that $B$ is trivial if there exists a 
$k$-algebra homomorphism $s:A\to B$, which is the left inverse to the residue 
homomorphism $B\to A$. Indeed, then $s\otimes 1:A\otimes _kk_n\to B$ is an 
isomorphism of $k_n$-algebras. 

Note that the quotient homomorphism $B\to B/t^nB$ defines the map 
$\de ^n(A)\to \de ^{n-1}(A)$. Denote by $\de _0^n(A)\subset \de ^n(A)$ the 
preimage  in $\de ^n(A)$ of the trivial deformation in $\de ^{n-1}(A)$. 

\begin{lemma} There exists a natural identification $\de ^n_0(A)=\de (A)$. 
In particular,  $\de ^n_0(A)$ has a natural structure of a $k$-vector space. 
\end{lemma}

\begin{pf} Let $B\in \de ^n(A)$ be such that $B/t^nB=A\otimes _kk_{n-1}$. 
Consider the obvious $k$-algebra homomorphism $A\to A\otimes _kk_{n-1}$ and 
the induced pullback diagram
$$
\begin{array}{ccrcccccc}
0 & \to & t^nB & \to & B^\prime & \to & A & \to & 0 \\
  &     & id \downarrow & & \downarrow & & \downarrow & & \\
0 & \to & t^nB & \to & B        & \to & A\otimes _kk_{n-1} & \to & 0
\end{array}
$$
Then $B^\prime $ represents an element in $\de (A)$. We get a map 
$\de _0^n(A)\to \de (A)$. 

The inverse map $\de (A) \to \de _0^n(A)$ is defined as follows. 
Given $B^\prime \in \de (A)$ consider the projection $A\otimes _kk_{n-1}
\to A$ and the corresponding pullback diagram
$$
\begin{array}{ccrcccccc}
0 & \to & A & \to & B & \to & A\otimes _kk_{n-1} & \to & 0 \\
  &     & id \downarrow & & \downarrow & & \downarrow & & \\
0 & \to & A & \to & B^\prime & \stackrel{p}{\to} & A & \to & 0
\end{array}
$$
Then $B$ is a $k_n$-algebra as follows:
$$t:(b^\prime ,0)\to (0,tp(b^\prime)),\quad t:(0,t^{n-1}a)\to 
(tp^{-1}(a),0).$$ 
This proves the lemma. 
\end{pf}

\begin{cor}
Assume that $\de (A)=0$. Then $\de ^n(A)=0$ for all $n$.
\end{cor}

\begin{pf} Induction on $n$ using the previous lemma.
\end{pf}

\subsection{Hochschild cohomology} The space $exal(A,M)$ has a well 
known cohomological description. Namely, there is a natural isomorphism 
$$exal(A,M)=\Ext ^2_{A^e}(A,M).$$
Let us recall how this isomorphism is defined. Consider the bar resolution 
$$...\stackrel{\partial _2}{\rightarrow}B_1\stackrel{\partial _1}
{\rightarrow}B_0
\stackrel{\partial _0}{\rightarrow}A\rightarrow 0,$$
where $B_i=A^{\otimes i+2}$ and 
$$\partial_i(a_0\otimes ...\otimes a_{i+1})=\sum_j
(-1)^ja_0\otimes ...\otimes a_ja_{j+1}\otimes ...a_{i+1}.$$
Note that $B_i$'s are naturally $A^e$-modules and the differentials 
$\partial _i$ are $A^e$-linear. Hence $B_\bullet \to A$ is a free resolution 
of the $A^e$-module $A$. Thus for any $A^e$-module $M$ 
$$H^\bullet \Hom_{A^e}(B_\bullet ,M)=\Ext^\bullet_{A^e}(A,M).$$
Note that $\Hom _{A^e}(B_i ,M)=\Hom _k(A^{\otimes i},M)$. 

Given an algebra extension 
$$0\to M\to B\to A\to 0$$
choose a $k$-linear splitting $s:A\to B$ and define a 2-cocycle 
$Z_s\in \Hom _k(A^{\otimes 2},M)$ by 
$$Z_s(a,b)=s(ab)-s(a)s(b).$$
Different $k$-splittings define cohomologous cocycles, hence we obtain a 
map $exal(A,M)\to \Ext^2_{A^e}(A,M)$ which is, in fact, an isomorphism.

The spaces $\Ext ^\bullet_{A^e}(A,M)$ are called the Hochschild cohomology 
groups of $A$ with coefficients in $M$. In particular, $\Ext _{A^e}^\bullet
(A,A)=HH^\bullet (A)$ is the usual Hochshild cohomology of $A$. Note that 
the space $\Ext ^0_{A^e}(A,M)=\Hom _{A^e}(A,M)$ coincides with the center 
$Z(M)$ of $M$:
$$Z(M)=\{ m\in M\vert am=ma\ \ \forall a\in A\} .$$
The space $\Ext ^1_{A^e}(A,M)$ classifies the outer derivations of $A$ 
into $M$. Namely, a map $d:A\to M$ is a {\it derivation} 
if $d(ab)=ad(b)+d(a)b$. 
It is called an inner derivation (defined by $m\in M$) if $d(a)=[a,m]$. 
Denote by 
$Der (A,M)$ (resp. $Inder(A,M)$) the space of derivations 
(resp. inner derivations). Then 
$$\Ext ^1_{A^e}(A,M)=Outder(A,M):=Der(A,M)/Inder(A,M).$$

\begin{remark} Consider the split extension $B=A\oplus M\in exal (A)$, 
i.e. the 
multiplication in $B$ is $(a,m)(a^\prime, m^\prime)=(aa^\prime, am^\prime+
ma^\prime)$. Then an automorphism of this extension is an algebra automorphism 
$\alpha \in Aut(B)$ of the form
$$\alpha (a,m)=(a,m+d(a)),$$
where $d:A\to M$ is a derivation. In other words the automorphism group 
of the trivial extension is the group $Der(A,M)$.
\end{remark}

\subsection{Deformation of sheaves of algebras} 
Let $X$ be a topological space and $\cA $ be a sheaf of $k$-algebras on $X$. 
Let $\cA ^o$ denote the sheaf of opposite $k$-algebras and 
$\cA ^e=\cA \otimes _k\cA ^o$. Given an $\cA ^e$-module $\cM $ we may repeat 
the above definition for algebras and modules to define the space 
of algebra extensions $exal(\cA,\cM)$. In particular, an algebra 
extension of $\cA $ by $\cM$ is represented by an exact sequence of sheaves of 
$k$-vector spaces 
$$0\to \cM \to \cB \stackrel{\epsilon}{\to} \cA \to 0$$
such that $\cB $ is a sheaf of $k$-algebras and $\epsilon$ is a homomorphism 
of sheaves of algebras  satisfying the properties of the Definition 
2.1 above. A split extension is the one admitting a homomorphism of sheaves 
of algebras $s:\cA \to \cB$ such that $\epsilon \cdot s=id$. In particular, 
a split extension must be split as an extension of sheaves of $k$-vector 
spaces. 

In case $\cM =\cA$ we may again define the set $\de ^n(\cA)$ of $n$-th order 
deformations of $\cA$, so that $\de ^1(\cA )=\de (\cA)=exal(\cA ,\cA)$. Let 
again $\de ^n_0(\cA )\subset \de ^n(\cA )$ be the subset consisting 
of $n$-th order deformations which are trivial up to order $n-1$. Then 
repeating the proof of Lemma 2.2 we get $\de ^n_0(\cA )=\de (\cA)$. In 
particular, $\de ^n_0(\cA )$ is naturally a $k$-vector space and $\de (\cA)=0$ 
implies $\de ^n(\cA)=0$ for all $n$.

\section{Review of Gerstenhaber-Schack construction}

In the paper [GS] the authors develop a deformation theory of presheaves of algebras 
on small categories. We will review their construction in a special case which is relevant
 to us. Namely let $X$ be a topological space and $\cU$ be the category of 
all (or some) open subsets of $X$. Let 
$\cA=\cA_{\cU}$ be a presheaf of algebras on $\cU$, i.e. $\cA$ is a contravariant functor 
from $\cU$ to the category of algebras. We denote by $k_{\cU}$ the constant presheaf 
of algebras: $k_{\cU}(U)=k$ for all $U\in \cU$.
Let $\cA -mod$ be the abelian category (of presheaves) of left $\cA$-modules. 
The presheaf of algebras $\cA ^e=\cA \otimes \cA ^o$ is defined in the obvious way: 
 $\cA ^e(U)=\cA(U)\otimes _k\cA ^0(U)$. 
In case $\cA=k_{\cU}$ for $\cM\in k_{\cU}-mod$ we denote 
$\Ext^i_{k_{\cU}}(k_{\cU},\cM )=H^i(\cU ,\cM ).$

Fix an $\cA$-bimodule $\cM$ (i.e. $\cM \in \cA^e-mod$). The group $exal(\cA ,\cM )$ 
is defined exactly as above in the case of a single algebra and its bimodule. We are going 
to give a natural description of the group $exal(\cA, \cM)$ in terms of homological algebra 
in the category of presheaves on $\cU$. In patricular, we will construct a canonical map 
$$exal(\cA ,\cM)\to \Ext ^2_{\cA ^e}(\cA ,\cM ).$$

First recall some constructions from [GS].

\subsection{Categorical simplicial resolution}

Let $\cC=\cC _{\cU}$ be a presheaf of algebras on $\cU$. Given $U\in \cU$ denote its 
inclusion $i_U:\{U\}\hookrightarrow \cU$. The obvious (exact) 
restriction functor 
$$i^*:\cC -mod \longrightarrow \cC (U)-mod,\quad \cK \mapsto \cK (U)$$
has a right exact left adjoint functor $i_{U!}:\cC(U)-mod\rightarrow \cC-mod$
$$i_{U!}K(V)=\begin{cases}
                \cC(V)\otimes _{\cC(U)}K, & \text{ if $V\subset U$},\\
                  0, & \text{ otherwise.}
\end{cases}$$
Thus if $K$ is a projective $\cC(U)$-module, then $i_{U!}K$ is a projective object in 
$\cC-mod$. In particular, the category $\cC-mod$ has enough projectives (it also 
has enough injectives (see [GS])). 

If the category $\cU$ has a final object $U$, then 
$\cC=i_{U!}\cC(U)$ is projective in $\cC-mod$. In patricular, then 
$$\Ext _{\cC}^i(\cC ,\cK )=0, \ \ \  \text{for all $\cK \in \cC-mod$, 
$i>0$.}$$

For $\cN\in \cC-mod$ define 
$$S(\cN):=\bigoplus_{U\in \cU}i_{U!}i_U^*\cN$$
with the canonical map
$$\epsilon _{\cN}:S(\cN)\rightarrow \cN.$$
Clearly $S$ is an endo-functor $S:\cC -mod\longrightarrow \cC -mod$ 
with a morphism of functors 
$\epsilon :S\to Id$. 

Define a diagram of functors 
$$...s_2\stackrel{\partial _1}{\rightarrow}s_1\stackrel{\partial _0}{\rightarrow}s_0
\stackrel{\partial _{-1}=\epsilon}{\longrightarrow}Id\to 0,$$
where $s_i=S^{i+1}$ and $\partial _i=\epsilon _{s_i}
 -S(\partial _{i-1})$. 
This diagram is a complex, i.e. $\partial _i\partial _{i-1}=0$, which is exact. 
So for $\cN \in \cC-mod$ we obtain a resolution
$$...\to s_1(\cN)\to s_0(\cN )\to \cN \to 0.$$
Explicitly we have 
$$s_k(\cN)=\bigoplus_{U_k\subset  ... \subset U_0}i_{U_k!}i_{U_k}^*...i_{U_0!}i_{U_0}^*\cN.$$
If $\cN$ is locally projective (i.e. $\cN(U)$ is a projective $\cC(U)-module$ for all 
$U\in \cU$), then the complex $s_\bullet(\cN)$ consists of projective objects in $\cC-mod$. 
So in this case for $\cM \in \cC-mod$ we have 
$$\Hom _{\cC}(s_\bullet (\cN),\cM)=\bbR \Hom_{\cC}^\bullet(\cN ,\cM ).$$

\subsection{Simplicial bar resolution}

Consider the bar resolution of the presheaf of algebras $\cA$:
$$...\to \cB _1\to \cB _0\to \cA,$$
where $\cB_i=\cA ^{\otimes i+2}$ (this is a direct analogue of the usual 
bar resolution for algebras described above).
 The presheaves $\cB _i$ are {\it locally} free $\cA ^e$-
modules, but usually not projective objects in $\cA ^e-mod$. So the simplicial 
resolution $s_\bullet \cB _\bullet$ of $\cB_\bullet$ is a double complex consisting of 
projective objects in $\cA ^e-mod$. For an $\cA ^e$-module $\cM $ denote by $T^{\bullet
\bullet}(\cM)$ the double complex $\Hom _{\cA ^e}(s_\bullet\cB_\bullet, \cM )$, and let 
$T^\bullet (\cM)=\Tot (T^{\bullet\bullet}(\cM ))$ be its total complex. We have 
$$\Ext ^i_{\cA ^e}(\cA ,\cM )=H^i(T^\bullet(\cM)).$$

Consider the double complex $T^{\bullet \bullet}(\cM)$. It looks like

\medskip

$$\begin{array}{lclc}
\uparrow &  & \uparrow &  \\
\prod\limits_U\Hom _k (\cA (U)\otimes \cA (U), \cM (U)) & \to & \prod
\limits_{V\subset U}
\Hom _k(\cA (U)\otimes \cA (U),\cM(V)) & \to  \\
\uparrow & & \uparrow &  \\
\prod\limits_{U}\Hom _k(\cA (U), \cM (U)) & \to & \prod\limits_{V\subset U}
\Hom _k(\cA (U),\cM (V)) & \to \\
\uparrow & & \uparrow &  \\
\prod\limits_U\Hom _k(k, \cM (U)) & \to & \prod\limits _{V\subset U}\Hom _k(k, \cM (V)) & \to ,
\end{array}$$

\medskip

\noindent where the left lower corner has bidegree $(0,0)$. 
The vertical arrows are the 
Hochshild differentials while the horizontal ones come 
from the simplicial resolution. 

Let $T_a^{\bullet \bullet}(\cM )\subset T^{\bullet \bullet}(\cM )$ be the sub- double 
complex which is the complement of the bottom row. 
Put 
$$T^\bullet_a(\cM)=Tot(T^{\bullet\bullet}_a(\cM));\quad H_a^n(\cA ,\cM):=
H^n(T_a^\bullet(\cM )).$$ 
Note that the complex $T^\bullet(\cM)/T_a^\bullet(\cM))$ is just 
$\Hom _k(s_\bullet (k_{\cU}),\cM)$. Hence we obtain the long exact sequence 
$$\begin{array}{cccccc}
\to & H_a^n(\cA ,\cM ) & \to & \Ext^n_{\cA ^e}(\cA ,\cM ) & \to & H^n(\cU ,\cM ) \\
\to & H_a^{n+1}(\cA ,\cM ) & \to & ... & & 
\end{array}
$$
In case $\cM $ is a symmetric $\cA$-bimodule, i.e. $am=ma$ for all $a\in \cA$, $m\in \cM$, 
the above sequence splits into short exact sequences ([GS],21.3)
$$0 \to H_a^n(\cA ,\cM )\to \Ext ^n_{\cA ^e}(\cA ,\cM) \to H^n(\cU ,\cM )\to 0.$$

\subsection{The isomorphism $exal(\cA ,\cM)\simeq H_a^2(\cA ,\cM)$}

Let the extension 
$$0\to \cM \to \cB \to \cA \to 0$$
represent an element in $exal(\cA ,\cM )$. Choose local $k$-linear splittings 
$s_U:\cA (U)\to \cB (U)$. Let us construct a 2-cocycle in 
$T^{\bullet\bullet}_a(\cM )$. 
Namely, put 
$$Z^{0,2}(a,b)=s_U(ab)-s_U(a)s_U(b), \quad U\in \cU,\ a,b\in \cA(U),$$
$$Z^{1,1}(a)=s_Vr^{\cA}_{U,V}(a)-r^{\cB}_{U,V}s_U(a), \quad V\subset U,\ a\in \cA (U),$$
where $r^{\cA}_{U,V}:\cA(U)\to \cA(V)$, $r^{\cB}_{U,V}:\cB (U)\to \cB (V)$ are the 
structure restriction maps of the presheaves $\cA $ and $\cB$. Then $(Z^{0,2},Z^{1,1})$ 
is a 2-cocycle in $T_a^{\bullet\bullet}(\cM)$ and the induced map 
$$exal(\cA ,\cM )\to H_a^2(\cA ,\cM)$$
is an isomorphism ([GS],21.4). 
The inverse isomorphism is constructed as follows. Let $(Z^{0,2},Z^{1,1})$ 
be a 2-cocycle in $T^{\bullet \bullet}_a(\cM)$. For each $U\in \cU$ put 
$\cB(U)=\cA(U)\oplus \cM(U)$ as a $k$-vector space; define the 
multiplication in $\cB(U)$ by $(a,m)(a^\prime ,m^\prime )
=(aa^\prime ,am^\prime + ma^\prime + Z^{0,2}(a,a^\prime ))$. We make $\cB$ 
the presheaf of algebras by defining the restriction maps $r_{U,V}^{\cB}
:\cB(U)\to \cB(V)$ to be $r_{U,V}^{\cB}(a,m)=(r^{\cA}_{U,V}(a),
r_{U,V}^{\cM}(m)+Z^{1,1}(a)).$

In particular, we obtain the 5-term exact 
sequence
$$\begin{array}{cccccc}
  & ...  & \to & \Ext^1_{\cA ^e}(\cA ,\cM ) & \to & H^n(\cU ,\cM ) \\
\to & exal(\cA ,\cM ) & \to & \Ext ^2_{\cA^e}(\cA,\cM) & \to & 
H^2(\cU, \cM) 
\end{array}
$$

\section{Admissible quasicoherent sheaves of algebras and bimodules.}

\begin{defn} Let $Z$ be a scheme and $\cA_Z$ be a sheaf of unital 
$k$-algebras on $Z$. 
We say that $\cA_Z$ is a {\it quasicoherent} sheaf of algebras if 
there is given a homomorphism of sheaves of unital $k$-algebras 
$\cO_Z\to \cA_Z$ which makes $\cA_Z$ a quasicoherent left $\cO_Z$-module. 
Note that $\cA_Z^o$ is then a quasicoherent right $\cO_Z$-module. 
Denote by $\mu(\cA_Z)\subset \cA_Z-mod$ the full subcategory of 
left $\cA_Z$-modules consisting of quasicoherent $\cO_Z$-modules
\end{defn}

Fix  a quasiprojective scheme $X$ over $k$ with a sheaf of unital 
$k$-algebras on $\cA_X$. Let $\cA_X^o$ be the sheaf of opposite algebras and 
$\cA^e_X=\cA_X\otimes _k\cA_X^o$. 
An $\cA_X$-module means a left $\cA_X$-module; 
an $\cA_X$-bimodule means an $\cA_X^e$-module. 
Put $Y=X\times X$ with the two projections $p_1,p_2:Y\to X$. We have the 
sheaves of algebras $p_1^{-1}\cA_X$ and $p_2^{-1}\cA_X^o$ on $Y$ and hence 
also their tensor product $p_1^{-1}\cA_X\otimes _kp_2^{-1}\cA_X^o$. 

Assume 
that $\cA_X$ is quasicoherent. Then we can take the quasicoherent inverse 
images $p_1^*\cA_X$ and $p_2^*\cA_X^o$ (using left and right 
$\cO_X$-structures respectively). Put 
$$\cA_Y^e:=p_1^*\cA_X\otimes _{\cO_Y}p_2^*\cA_X^o.$$
Note that for affine open $U,V\subset X$, $\cA_Y^e(U\times V)=
\cA_X(U)\otimes _k\cA_X(V)$. This is a quasicoherent sheaf on $Y$ with 
a natural morphism of quasicoherent sheaves 
$$\beta :\cO_Y\to \cA_Y^e,$$
which sends $1$ to $1\otimes 1$. We also have the obvious morphism of 
sheaves of $k$-vector spaces
$$\gamma :p_1^{-1}\cA_X\otimes _kp_2^{-1}\cA_X^o\to \cA_Y^e.$$

\begin{defn} We say that the quasicoherent sheaf of algebras $\cA_X$ 
satisfies condition (*) if $\cA_Y^e$ has a structure of a sheaf of algebras 
so that $\beta $ and $\gamma $ are morphisms of sheaves of algebras. 
\end{defn}

Note that if $\cA_X$ satisfies condition (*) then, in particular, 
 $\cA_Y^e$ is a quasicoherent 
sheaf of algebras on $Y$. It seems that the algebra structure on $\cA_Y^e$ 
as required in the condition (*), if it exists, should be unique. In any case, 
there is a canonical such structure in all examples that we have in mind. 

\bigskip

\noindent{\bf Examples.} 1. The condition (*) holds if the sheaf of algebras 
$\cA_X$ is commutative. More generally, if the image of $\cO_X$ lies in the 
center of $\cA_X$. 

2. Assume that $char(k)=0$ and $X$ is smooth. Then (*) holds for the sheaf 
$\cA_X=D_X$ of differential operators on $X$. In this case 
$$p_1^*D_X\otimes _{\cO_Y}p_2^*D_X=D_Y.$$
Let $\omega _X$ be the dualizing sheaf on $X$. Then $D_X^o=\omega _X
\otimes _{\cO_X}D_X\otimes _{\cO_X}\omega _X^{-1}$ and hence 
$$\cA^e_Y=p_1^*D_X\otimes _{\cO_Y}p_2^*D_X^0=p_2^*\omega _X\otimes _{\cO_Y} 
D_Y\otimes _{\cO_Y}p_2^*\omega _X^{-1}.$$

\bigskip

Let $\cM_X$ be an $\cA_X$-bimodule. Then, in particular, $\cM_X$ is an 
$\cO_X$-bimodule. 

\begin{defn} We say that $\cM_X$ satisfies the condition $(\star)$ if for an 
open affine $U\subset X$ and $f\in \cO(U)$ we have 
$$\cM_X(U_f)=\cO(U_f)\otimes _{\cO(U)}\cM_X(U)\otimes _{\cO(U)}
\cO(U_f).$$
\end{defn}

\begin{remark} The sheaves of algebras $\cA_X$ in Examples 1,2 above 
satisfy the condition $(\star)$ when considered as $\cA_X$-bimodules.
\end{remark}

\begin{lemma} Let $\cA_X$ be a quasicoherent sheaf of algebras which 
satisfies the condition (*), and let $\cM_X$ be an $\cA_X$-bimodule which 
satisfies the condition $(\star )$. Then $\cM_X$ defines a (unique up to an 
isomorphism) $\cA^e_Y$-module $\tilde{\cM}_Y$ on $Y$ such that for an affine 
open $U\subset X$  
$$\tilde{\cM}_Y(U\times U)=\cM_X(U).$$
We have $\tilde{\cM}_Y\in \mu(\cA^e_Y)$. 
\end{lemma}

\begin{pf} Choose an affine open covering $\{U\}$ of $X$. Then the affine 
open subsets $U\times U$ form a covering of $Y$. Fix one such subset 
$V=U\times U$. The sheaf of algebras $\cA^e_Y$ is quasicoherent, hence by 
Serre's theorem below we have the equivalence of categories
$$\mu(\cA^e_V)\simeq \cA_Y^e(V)-mod.$$
The sheaf $\cM_X$ defines an $\cA^e_Y(V)=\cA_X(U)\otimes _k\cA_X(U)$-module 
$\cM_X(U)$, hence defines a quasicoherent $\cA_V^e$-module 
$\tilde{\cM}_V$. If $V^\prime =U^\prime\times U^\prime \subset V$, then 
the condition $(\star)$ for $\cM_X$ implies that $\tilde{\cM}_V
\vert_{V^\prime}=\tilde{\cM}_{V^\prime}$. Hence the local sheaves glue 
together into a global quasicoherent $\cA^e_Y$-module $\tilde{\cM}_Y$. 
The last assertion is obvious. 
\end{pf}

\begin{thm} Let $Z=SpecC$ be an affine scheme, $\cA_Z$ -- a quasicoherent 
sheaf 
of algebras on $Z$. Put 
$A=\Gamma (X,\cA_X)$. Then the functor of global sections $\Gamma $ is an 
equivalence of categories 
$$\Gamma :\mu(\cA_Z)\to A-mod.$$
Its inverse is $\Delta$ defined by 
$$\Delta (M)=\cA_Z\otimes _AM.$$
Both $\Gamma $ and $\Delta $ are exact functors. 
\end{thm}

\begin{pf} The point is that for an $A$-module $M$ the quasicoherent sheaf 
$\Delta (M)$ is indeed an $\cA_Z$-module. The rest follows easily from the 
classical Serre's theorem about the equivalence
$$qcoh(Z)\simeq C-mod.$$
\end{pf}

\begin{defn} We call a quasicoherent sheaf of algebras $\cA_X$ 
{\it admissible} if it satisfies conditions (*) and $(\star)$ 
(as a bimodule over itself). We call an $\cA_X$-bimodule 
$\cM_X$ 
{\it admissible} in it satisfies condition $(\star)$. We say that 
$(\cA_X,\cM_X)$ is an admissible pair if both $\cA_X$ and $\cM_X$ are 
admissible. 
\end{defn}

\begin{remark} The sheaf of algebras $\cA_X$ as in Examples 1,2 above is 
admissible.
\end{remark}

Let us summarize our discussion in the following corollary.

\begin{cor} Let $(\cA_X, \cM_X)$ be an admissible pair. Then

i) $\cA_X$ defines  is a quasicoherent sheaf of algebras $\cA_Y^e$ 
on $Y$ such that for 
affine open $U,V\subset X$, $\cA_Y^e(U\times V)=\cA_X(U)\otimes _k
\cA_X(V)^o$; 

ii) $\cM_X$ defines a sheaf $\tilde{\cM}_Y\in \mu(\cA_Y^e)$ such 
that for affine open $U\subset X$, $\tilde{\cM}_Y(U\times U)=\cM_X(U).$
\end{cor}

\begin{pf} This follows immediately from Definition 4.2 and Lemma 4.5.
\end{pf}

We will be able to give a cohomological interpretation of the group 
$exal(\cA_X,\cM_X)$ for an admissible pair $(\cA_X,\cM_X)$.

\section{Cohomological description of the group $exal(\cA_X,\cM_X)$ for 
an admissible pair $(\cA_X,\cM_X)$.}

Let $X$ be a quasiprojective scheme over $k$ and $(\cA_X,\cM_X)$ be an 
admissible pair. We will consider the group $exal(\cA_X,\cM_X)$ of 
algebra extensions of $\cA_X$ by $\cM_X$. Note that if an exact sequence 
$$0\to \cM_X \to \cB _X \to \cA _X \to 0$$
is such an extension, then we do not require the sheaf $\cB_X$ to be 
quasicoherent, or even an $\cO_X$-module. 

Denote by  $\cU =Aff(X)$ be the category of all affine 
open subsets of $X$. Given a sheaf 
$\cF _X$ on $X$ we denote by $j_X^*\cF _X$ the presheaf on $\cU$, 
which is obtained 
by restriction of  $\cF _X$ to affine open subsets. We will usually denote 
$j^*_X\cF _X=\cF _{\cU}$ if it causes no confusion. 
In particular, we obtain presheaves 
of algebras $\cA _{\cU}=j^*_X\cA _X$, $\cA ^e_{\cU}:=\cA _{\cU}\otimes \cA _{\cU}^o$ $(\cA ^e_{\cU}\neq j^*_X\cA ^e_X)$. 

\begin{lemma} Then there is a  natural map 
$exal(\cA _X, \cM _X)\rightarrow exal(\cA _{\cU},\cM _{\cU})$ which is an isomorphism. 
In particular, $\de (\cA _X)=\de (\cA _{\cU})$. 
\end{lemma}

\begin{pf} Given an exact sequence of sheaves on $X$
$$0\to \cM_X \to \cB _X \to \cA _X \to 0,$$
which represents an element in $exal(\cA _X,\cM _X)$ we obtain the corresponding sequence 
$$0\to \cM _{\cU}\to \cB _{\cU} \to \cA_{\cU}\to 0$$
of presheaves on $\cU$. This last sequence is exact because $\cM _X$ is 
quasi-coherent. 
Hence it represents an element in $exal(\cA _{\cU},\cM _{\cU})$. So we obtain a map 
$$exal(\cA _X,\cM _X)\rightarrow exal(\cA _{\cU},\cM _{\cU}).$$

Vice versa, let 
$$0 \to \cM _{\cU }\to \cB _1 \to \cA _{\cU} \to 0$$
represent an element in $exal(\cA _{\cU}, \cM_{\cU})$. Denote by ${}^+$ the 
(exact) functor which associates to a presheaf on $\cU$ the corresponding 
sheaf on $X$. Then $(\cA _{\cU})^+=\cA _X$, $(\cM _{\cU})^+=\cM _X$ and hence 
we obtain an exact sequence 
$$0 \to \cM _X \to \cB_1^+ \to \cA _X \to 0$$
which defines an element in $exal(\cA _X,\cM_X)$. This defines 
the inverse map 
$$exal(\cA _{\cU},\cM_{\cU})\rightarrow exal(\cA _X,\cM _X).$$
\end{pf}

Let $D^b(\cA_Y^e)$ and $D^b(\cA^e_{\cU})$ denote the bounded 
derived categories of $\cA_Y^e-mod$ and $\cA^e_{\cU}-mod$ respectively. 
Let $D^b_{\mu(\cA^e_Y)}(\cA_Y^e)\subset D^b(\cA_Y^e)$ be the full 
subcategory consisting of complexes with cohomologies in 
$\mu(\cA^e_Y)$. 
Denote by $j^*_Y: \cA^e_Y-mod \longrightarrow \cA^e_{\cU}-mod$ the left exact 
functor defined by $j^*_Y(\cF)(U):=\cF(U\times U)$, $U\in \cU$. Consider 
its derived functor 
$$\bbR j^*_Y:D^b(\cA^e_Y)\longrightarrow D^b(\cA^e_{\cU}).$$

\begin{thm} The functor 
$$\bbR j^*_Y:D^b_{\mu(\cA^e_Y)}(\cA^e_Y)\longrightarrow 
D^b(\cA^e_{\cU})$$
is fully faithful. Equivalently, for $\cM,\cN\in \mu (\cA _Y^e)$ the map 
$$j^*_Y:\Ext ^n_{\cA ^e_Y}(\cM ,\cN)\rightarrow \Ext ^n_{\cA^e_{\cU}}
(j^*_Y\cM,j^*_Y\cN)$$
is an isomorphism for all $n$. 
\end{thm}

\begin{prop} The map 
$$j^*_X:H^n(X,\cM_X)\rightarrow H^n(\cU, \cM_{\cU})$$   
is an isomorphism for all $n$. 
\end{prop}

Let us first formulate some immediate corollaries of the theorem and the 
proposition.

\begin{cor} There 
exists a natural exact sequence
$$\Ext ^1_{\cA ^e_Y}(\tilde{\cA}_Y,\tilde{\cM}_Y)\to H^1(X,\cM_X)\to 
exal(\cA _X,\cM_X )
\to \Ext ^2_{\cA ^e_Y}
(\tilde{\cA}_Y,\tilde{\cM}_Y)\to H^2(X,\cM_X).$$
In particular, if $X$ is affine then 
$exal(\cA _X,\cM_X )=\Ext^2_{\cA ^e_Y}(\tilde{\cA}_Y,\tilde{\cM}_Y)$. 
If $\cM_X$ is a symmetric $\cA_X$-bimodule, then we get a short exact sequence 
$$0\to exal(\cA _X,\cM _X)\to \Ext ^2_{\cA^e_Y}(\tilde{\cA}_Y,\tilde{\cM}_Y)
\to H^2(X,\cM_X)\to 0.$$
\end{cor}

\begin{pf} Indeed, this follows from Lemma 5.1, Theorem 5.2, Proposition 5.3 
and results of Section 3. 
\end{pf}

Recall the following theorem of J.~Bernstein.

\begin{thm}([Bo]) Let $Z$ be a quasicompact separated scheme, $\cC_Z$ -- 
a quasicoherent sheaf of algebras 
on $Z$. Then the 
natural functor 
$$\theta:D^b(\mu(\cC_Z))\to D^b_{\mu(\cC_Z)}(\cC_Z)$$ 
is an equivalence of 
categories.
\end{thm}   

\begin{cor} Assume that $X$ is affine. Then 
$$exal(\cA_X,\cM_X)\simeq exal(\cA_X(X),\cM_X(X)).$$ 
\end{cor}

\begin{pf} Put $\cA_X(X)=A$, $\cM_X(X)=M$. We have 
$$exal(A,M)=\Ext ^2_{A\otimes A^o}(A,M).$$
By Serre's theorem 
$$\Ext ^2_{A\otimes A^o}(A,M)=\Ext ^2_{\mu(\cA^e_Y)}(\tilde{\cA}_Y,
\tilde{\cM}_Y).$$
By Bernstein's theorem
$$\Ext ^2_{\mu(\cA^e_Y)}(\tilde{\cA}_Y,
\tilde{\cM}_Y)=\Ext ^2_{\cA^e_Y}(\tilde{\cA}_Y,
\tilde{\cM}_Y).$$
Finally, by Corollary 5.4 above
$$\Ext ^2_{\cA^e_Y}(\tilde{\cA}_Y,
\tilde{\cM}_Y)=exal(\cA_X,\cM_X).$$
\end{pf}

\noindent{\it Question.} Under the assumptions of the last corollary let 
$\cB$ be a sheaf of algebras on $X$ representing an element in $exal(\cA_X, 
\cM_X)$. Is $\cB =\cA_X \oplus \cM_X$ as a sheaf of $k$-vector spaces?

\section{Proof of Theorem 5.2 and Proposition 5.3.}

\medskip

\noindent{\it Proof of Proposition 5.3.}
Let $k_{\cU}$ be the constant presheaf on $\cU$ and $s_\bullet (k_{\cU})\to 
k_{\cU}$ be its categorical simplicial resolution (Section 3). It is 
a projective resolution of $k_{\cU}$, which consists of direct sums of 
presheaves $i_{U!}k$. Hence 
$$H^i(\cU,\cM_{\cU})=\Ext ^i(k_{\cU},\cM_U)=H^i\Hom^\bullet 
(s_{\bullet}(k_{\cU}),
\cM_{\cU}).$$

Consider the exact functor $(\cdot )^+$ from the category of presheaves on 
$\cU$ to the category on sheaves on $X$. Then $k_{\cU}^+=k_X$ -- the constant 
sheaf on $X$. The functor $(\cdot )^+$ preserves direct sums and 
$(i_{U!}k)^+=k_U$ -- the extension by zero of the constant sheaf on $U$. 
Since $\cM_X$ is quasicoherent, for an affine open $U\subset X$ we have 
$H^i(U,\cM_X)=0$ for all $i>0$. 
Thus 
$$H^i(X,\cM_X)=\Ext ^i(k_X,\cM_X)=H^i\Hom ^\bullet(s_\bullet(k_{\cU})^+,
\cM_X).$$

It remains to notice that 
$$\Hom (k_U,\cM_X)=\Gamma (U,\cM_X)=\Hom (i_{U!}k,\cM_{\cU}).$$
This completes the proof of the proposition.

\medskip

\noindent{\it Proof of Theorem 4.2.} 

Let us formulate a general statement which will imply the theorem. 
Let $Z$ be a quasicompact separated scheme over $k$. Let $Aff(Z)$ be 
the category of affine open subsets of $Z$ and 
$\cW\subset 
Aff(Z)$ be a full subcategory which is closed under intersections and 
constitutes a covering of $Z$. Let $\cA_Z$ be a quasicoherent sheaf 
of algebras on $Z$. Denote by $\cA_W$ 
the corresponding presheaf of algebras on $\cW$. Let 
$$j_Z^*:\cA_Z-mod \longrightarrow \cA_W-mod$$ 
be the natural (left exact) 
restriction functor.

\begin{prop} In the above notation the derived functor 
$$\bbR j_Z^*:D^b_{\mu(\cA_Z)}(\cA_Z)\to D^b(\cA_W)$$
is fully faithful. 
\end{prop}

\begin{pf} By Bernstein's theorem the natural functor 
$$\theta :D^b(\mu(\cA_Z))\to D^b_{\mu(\cA_Z)}(\cA_Z)$$
is fully faithful. So it suffices to prove that the composition 
$\bbR j_Z^*\cdot \theta$ is fully faithful. The functor $j^*_Z:\mu(\cA_Z)
\to \cA_W-mod$ is exact. 
Let $\cM,\cN \in \mu(\cA_Z)$. It suffices to prove that the map 
$$j^*_Z:\Ext^\bullet_{\mu(\cA_Z)} (\cM,\cN)\to 
\Ext^\bullet (j^*_Z\cM,j^*_Z\cN)$$
is an isomorphism. 

\medskip
\noindent{\it Step 1.} Assume that $Z$ is affine and $Z\in \cW$. Then by 
Serre's theorem $\mu(\cA_Z)\simeq \cA_Z(Z)-mod$. Replacing $\cM$ by a 
left free resolution we may assume that $\cM=\cA_Z$. But then 
$$\Ext ^i(\cA_Z,\cN)=\Ext ^i(\cA_Z(Z),\cN(Z))=\begin{cases}
                                        \cN(Z), \text{if $i=0$}\\
                                         0, \text{otherwise}
\end{cases}$$
 
On the other hand $j^*_Z\cA_Z=\cA_{\cW}$ is a projective object in 
$\cA_{\cW}-mod$ (Section 3) and 
$$\Hom (\cA_{\cW},j_Z^*\cN)=\Hom (\cA_{\cW}(Z),j^*_Z\cN(Z))=\cN(Z).$$
So we are done.

\medskip
\noindent{\it Step 2. Reduction to the case when $Z$ is affine.}

Let $i_U:U\hookrightarrow Z$ be an embedding of some $U\in \cW$. 
Denote by $\cA_U$ the restriction $\cA_Z\vert _U$. We have two (exact) 
adjoint functors $i^*_U:\mu(\cA_Z)\to \mu(\cA_U)$, $i_{U*}:
\mu(\cA_U)\to \mu (\cA_Z)$. The functor $i_{U*}$ preserves injectives. 

Choose a finite covering $Z=\bigcup U_j$, $U_j\in \cW$. Then the natural map 
$$\cN \to \bigoplus_ji_{U_j*}i^*_{U_j}\cN$$
is a monomorphism. So we may assume that $\cN=i_{U*}\cN_U$ for some 
$U\in \cW$ and $\cN_U\in \mu(\cA_U)$. Then we have 
$$\Ext ^\bullet(\cM ,i_{U*}\cN_{U})=\Ext ^\bullet(i^*_U\cM,\cN_U).$$

We need a similar construction on the other end. Let 
$\tilde{i}_U:\cW_U\hookrightarrow \cW$ be the embedding of the full 
subcategory $\cW_U=\{ V\in \cW\vert V\subseteq U\}$. Let 
$\cA _{\cW_U}$ be the restriction of $\cA_{\cW}$ to $\cW_U$. We have 
the obvious functor $\tilde{i}_U^*:\cA_{\cW}-mod\longrightarrow 
\cA_{\cW_U}-mod$ 
and its right adjoint $\tilde{i}_{U*}$ defined by 
$$\tilde{i}_{U*}(\cK)(V):=\cK(V\cap U).$$
Both $\tilde{i}_U^*$ and $\tilde{i}_{U*}$ are exact and $i_{U*}$ 
preserves injectives. For $\cK \in \cA_{\cW_U}$, $\cL \in \cA_{\cW}$ 
we have 
$$\Ext^\bullet (\tilde{i}^*_U\cL,\cK)=\Ext ^\bullet(\cL,\tilde{i}_{U*}
\cK).$$

Note that the following diagrams commute
$$\begin{array}{ccc}
\mu(\cA_Z) & \stackrel{i^*_U}{\longrightarrow} & \mu(\cA_U) \\
j^*_X\downarrow & & \downarrow j^*_U \\
\cA_{\cW}-mod & \stackrel{\tilde{i}^*_U}{\longrightarrow} & \cA_{\cW_U}-mod 
\end{array}$$

$$\begin{array}{ccc}
\mu(\cA_Z) & \stackrel{i_{U*}}{\longleftarrow} & \mu(\cA_U) \\
j^*_X\downarrow & & \downarrow j^*_U \\
\cA_{\cW}-mod & \stackrel{\tilde{i}_{U*}}{\longleftarrow} & \cA_{\cW_U}-mod 
\end{array}$$
(here $j_U^*$ is the obvious restriction functor). Hence the following 
diagram commutes as well
$$\begin{array}{lcl}
\Ext^\bullet(\cM,\cN)= & \stackrel{j^*_Z}{\longrightarrow} & 
\Ext^\bullet(j^*_Z\cM,j_Z^*\cN) =\\
\Ext^\bullet(\cM,i_{U*}\cN_U)= & & \Ext^\bullet (j^*_Z\cM,
\tilde{i}_{U*}j^*_U\cN_U)=\\
\Ext ^\bullet(i^*_U\cM,\cN_U)= & \stackrel{j^*_U}{\longrightarrow} & 
\Ext ^\bullet(j^*_Ui^*_U\cM,j^*_U\cN_U).
\end{array}$$
But $j^*_U$ is an isomorphism by Step 1 above. Hence $j^*_Z$ is also an 
isomorphism. 
\end{pf}

\section{A spectral sequence}

Let $X$ be a quasiprojective variety and $(\cA_X,\cM_X)$ be an admissible 
pair. 
For $\cN_1,\cN_2\in \mu(\cA^e_Y)$ we will construct a spectral sequence which abuts to 
$\Ext^\bullet_{\cA^e_Y}(\cN_1,\cN_2)$. In particular we will get an insight 
into the group $\Ext ^2_{\cA^e_Y}(\tilde{\cA}_Y,\tilde{\cM}_Y)$. 

\begin{lemma} Any object in $\mu (\cA^e_Y)$ is a quotient of a locally free 
$\cA^e_Y$-module.
\end{lemma}

\begin{pf} Let $\cK\in \mu(\cA^e_Y)$. Consider $\cK$ as a quasi-coherent $\cO_Y$-module. 
As such it is a quotient of a locally free $\cO_Y$-module $Q$ (we can take 
$Q=\oplus\cO_Y(-j)$). Then the $\cA_Y$-module $\cA_Y\otimes _{\cO_Y}Q$ is locally free 
and surjects onto $\cK$.
\end{pf}

Let $P_\bullet \to \cN_1$ be a resolution of $\cN_1$ consisting of locally free $\cA_Y^e$
-modules. From the proof of the last lemma it follows that there exists an affine 
covering $\cV$ of $Y$ such that for each $V\in \cV$ and each $P_{-t}$ the restriction 
$P_{-t}\vert _V$ is a free $\cA^e_V$-module. 
We may (and will) assume that each $V\in \cV$ is of the form $U\times U$ for 
$U$ from an affine open covering $\cU$ of $X$. 
Choose one such affine covering $\cV$. 
Let $\check{C}_\bullet(P_\bullet)\rightarrow P_\bullet$ be the corresponding 
\v{C}ech resolution of $P_\bullet$. This is a double complex consisting of $\cA_Y^e$-modules, 
which are extensions by zero from affine open subsets $V$ of free $\cA^e_V$-modules. 
Thus 
$$H^\bullet_{\cA^e_Y} \Hom(Tot(\check{C}_\bullet(P_\bullet)),\cN_2)=
\Ext^\bullet_{\cA^e_Y}(\cN_1,\cN_2).$$

The natural filtration of the double complex $\check{C}_\bullet(P_\bullet)$ gives rise to the 
spectral sequence with the $E_2$-term 
$$E_2^{p,q}=\check{H}^p(\cV, \cE xt^q_{\cA^e_Y}(\cN_1,\cN_2)),$$
which abuts to $\Ext ^{p+q}_{\cA^e_Y}(\cN_1,\cN_2)$.  
 
In particular, in case $\cN_1=\tilde{\cA}^e_Y$, $\cN_2=\tilde{\cM}_Y$ this 
spectral sequence defines a filtration of the group 
$\Ext ^2_{\cA^e_Y}(\tilde{\cA}^e_Y,\tilde{\cM}_Y)$. Namely there are maps

$$\alpha _1:\Ext^2_{\cA_Y^e}(\tilde{\cA}_Y^e,\tilde{\cM}_Y^e)
\rightarrow \check{H}^0(\cV,
\cE xt^2_{\cA_Y^e}(\tilde{\cA}^e_Y,\tilde{\cM}_Y^e)),$$
$$\alpha _2:\ker(\alpha _1)\rightarrow \check{H}^1(\cV,\cE xt^1_{\cA_Y^e}
(\tilde{\cA}_Y^e,\tilde{\cM}_Y^e)),$$
$$\alpha _3:\ker(\alpha _2)\rightarrow \check{H}^2(\cV,\cE xt^0_{\cA_Y^e}
(\tilde{\cA}_Y^e,\tilde{\cM}_Y^e).$$

Recall that for $V=U\times U\in \cV$ by Bernstein's and Serre's theorems 
respectively we have 
$$\Gamma (V,\cE xt^q_{\cA_Y^e}(\tilde{\cA}_Y^e,\tilde{\cM}^e_Y))$$
$$=\Gamma (V,\cE xt^q_{\mu(\cA_Y^e)}(\tilde{\cA}_Y^e,\tilde{\cM}^e_Y))$$
$$=\Ext^q_{\cA_X(U)\otimes \cA_X^o(U)}(\cA_X(U),\cM_X(U)).$$

\subsection{Cohomological analysis of the group $exal(\cA_X,\cM_X)$}
Consider the exact sequence 
$$H^1(X,\cM_X)\stackrel{\epsilon}{\rightarrow} exal(\cA_X,\cM_X)
\stackrel{\rho }{\rightarrow} \Ext^2_{\cA_Y^e}(\tilde{\cA}^e_Y,
\tilde{\cM}_Y).$$
Let us describe the morphisms $\epsilon$ and $\rho$ explicitly.   

Since $\cM_X$ is quasi-coherent the cohomology group $H^1(X,\cM_X)$ 
is isomorphic to the 
\v{C}ech cohomology $\check{H}^1(\cU,\cM_X)$. 
Given a 1-cocycle $\{m_{ij}\in \cM_X(U_i\cap U_j)\vert U_i,U_j\in \cU\}$ 
define an algebra 
extension 
$$0\to \cM_X \to \cB \to \cA_X \to 0$$ as follows: on each $U\in \cU$ 
the sheaf $\cB\vert_U$ 
is a direct sum of sheaves $\cM_X\vert_U$ and $\cA_X\vert_U$ with the multiplication 
$$(m,a)(m^\prime,a^\prime)=(ma^\prime+am^\prime,aa^\prime).$$
That is, locally $\cB$ is a split extension. Define the glueing algebra automorphisms 
$$\phi _{ij}:\cB_{U_i\cap U_j}\stackrel{\sim}{\rightarrow} \cB_{U_i\cap U_j},\quad 
\phi _{ij}(m,a)=(m+[a,m_{ij}],a).$$
This defines the map $\epsilon: H^1(X,\cM_X)\to exal(\cA_X,\cM_X)$. 

Now assume that an algebra extension $\cB$ represents an element in $exal(\cA_X,\cM_X)$. 
Consider $\rho (\cB)\in \Ext^2_{\cA_Y^e}(\tilde{\cA}_Y,\tilde{\cM}_Y)$ 
and assume 
that $\alpha _1
(\rho (\cB))=0$, i.e. locally $\cB$ is a split extension. Thus for $U\in \cU$ 
we have 
$$\cB(U)=\cM_X(U)\oplus \cA_X(U)$$
with the multiplication
$$(m,a)(m^\prime,a^\prime)=(ma^\prime+am^\prime,aa^\prime)$$
and with the glueing given by algebra automorphisms 
$$\phi _{ij}:\cB(U_i\cap U_j)\stackrel{\sim}{\rightarrow}\cB(U_i\cap U_j), \quad 
\phi _{ij}(m,a)=(m+\delta _{ij}(a),a),$$
where $\delta _{ij}:\cA_X(U_i\cap U_j)\to \cM_X (U_i\cap U_j)$ 
is a derivation. 
For an affine open $U\subset X$ the space 
$$\Ext ^1_{\cA_X(U)\otimes \cA_X^o(U)}(\cA_X(U),\cM_X(U))$$ 
is the space of outer derivations 
$\cA_X(U)\to \cM_X(U)$. 
The collection $\{\delta _{ij}\}$ defines an element in 
$\check{H}^1(\cV,\cE xt^1_{\cA_Y^e}(\tilde{\cA}_Y,\tilde{\cM}_Y))$, 
which is equal to 
$\alpha _2(\rho (\cB))$. 

Assume now that $\alpha _2(\rho(\cB))=0$.
Then there exist elements $\delta _i\in \Ext ^1_{\cA_X(U_i)
\otimes \cA_X(U_i)^0}
(\cA_X(U_i),\cM_X(U_i))$ such that $\delta _{ij}=\delta _i-\delta _j$. 
Changing the local trivializations of $\cB$ by the  derivations 
$\delta _i$'s we may assume that $\delta _{ij}$'s are inner derivations. 
Choose 
$m_{ij}\in \cM_X(U_i\cap U_j)$ so that $\delta_{ij}(a)=[a,m_{ij}].$ The collection 
$\{m_{ij}\}$ defines a 1-cochain in $\check{C}(\cU,\cM_X)$. Its coboundary is a 
2-cocycle which consists of central elements $m_{ijk}\in \cM_X(U_i\cap U_j\cap U_k)$. Thus 
it defines an element in $\check{H}^2(\cV, \cH om_{\cA_Y^e}
(\tilde{\cA}_Y,\tilde{\cM}_Y))$. It is 
equal to $\alpha _3(\rho (\cB))$. 

\section{Examples} 

Let $X$ be a smooth complex quasiprojective variety. Let 
$\delta :X\hookrightarrow Y=X\times X$ be the diagonal embedding, 
$\Delta=\delta (X)$ -- the diagonal, and $p_1,p_2:Y\to X$ be the two 
projections.

\subsection{Deformation of the structure sheaf} Let $\cA_X=\cM_X=\cO_X$. 
Then $\cA^e_Y=\cO_Y$, $\tilde{\cA}_Y=\delta _*\cO_X$. 
Since the $\cO_X$-bimodule 
$\cO_X$ is symmetric we have the short exact sequence 
$$0\to \de (\cO_X)\to \Ext ^2_{\cO_Y}(\delta _*\cO_X,\delta _*\cO_X)\to 
H^2(X,\cO_X)\to 0.$$
Assume that $X$ is projective. 
By the Hodge decomposition ([GS],[S])
$$\Ext ^2_{\cO_Y}(\delta _*\cO_X,\delta _*\cO_X)=H^0(X,\wedge ^2T_X)
\oplus H^1(X,T_X)\oplus H^2(X,\cO_X).$$
The above short exact sequence identifies $\de (\cO_X)$ with 
$H^0(X,\wedge ^2T_X)
\oplus H^1(X,T_X)$. The summand $H^1(X,T_X)$ corresponds to 
the first order deformations of the variety $X$ by Kodaira-Spencer theory, 
i.e. to ``commutative'' deformations of $\cO_X$, while the summand 
$H^0(X,\wedge ^2T_X)$ corresponds to ``noncommutative'' deformations. 

\subsection{Deformations of the sheaf of differential operators} 
Let $\cA_X=\cM_X=D_X$ -- the sheaf 
of (algebraic) differential operators on $X$. 
Let $\omega _X$ be the dualizing sheaf on $X$. Then 
$$D_X^o=\omega _X\otimes _{\cO_X}D_X\otimes _{\cO_X}\omega^{-1}_X.$$
We have $D_Y=p^*_1D_X\otimes _{\cO_Y}p^*_2D_X$,
and hence
$$D_Y^e=p_1^*\omega _X\otimes _{\cO_Y}D_Y
\otimes _{\cO_Y}p_2^*\omega^{-1}_X.$$
The functor $\tau :M\mapsto p_1^*\omega _X\otimes _{\cO_Y}M$ is an 
equivalence of categories 
$$\tau :D_Y-mod \longrightarrow D_Y^e-mod.$$
Denote by $\delta _+:D_X-mod \longrightarrow D_Y-mod$ the functor of direct 
image ([Bo]). Then 
$$\tilde{D}_Y=\tau (\delta _+\cO_X).$$

Let $X^{an}$ denote the variety $X$ with the 
classical topology.

\begin{prop} There is a natural isomorphism
$$\Ext^\bullet_{D^e_Y}(\tilde{D}_Y,\tilde{D}_Y)\simeq H^\bullet(X^{an},\bbC).$$
\end{prop}

\begin{pf} By the above remarks
$$\Ext^\bullet_{D^e_Y}(\tilde{D}_Y,\tilde{D}_Y)=
\Ext^\bullet_{D_Y}(\delta _+\cO_X,\delta _+\cO_X).$$

Let $D_{\Delta}^b(D_Y)$ be the full subcategory of 
$D^b(D_Y)$ consisting of complexes with cohomologies supported 
on $\Delta$. By Kashiwara's theorem the direct image functor 
$$\delta _+:D^b(D_X)\longrightarrow D^b_{\Delta}(D_Y)$$
is an equivalence of categories (see [Bo]). Thus, in particular,
$$\Ext^\bullet _{D_X}(\cO_X,\cO_X)\simeq 
\Ext^\bullet _{D_Y}(\delta _+\cO_X,\delta _+\cO_X).$$
 On the 
other hand by (a special case of) the Riemann-Hilbert correspondence 
$$\Ext^\bullet _{D_X}(\cO_X,\cO_X)\simeq H^\bullet(X^{an},\bbC).$$
\end{pf}

\begin{cor} Let $X$ be a smooth complex quasi-projective variety. Then 
we have an exact sequence
$$H^1(X^{an},\bbC)\to H^1(X,D_X)\to \de (D_X)\to H^2(X^{an},\bbC)
\to H^2(X,D_X).$$
If $X$ is $D$-affine (for example $X$ is affine) then 
$$\de (D_X)=H^2(X^{an},\bbC).$$
\end{cor}

\begin{pf} The first part follows immediately from Proposition 8.1 and 
Corollary 5.4.  
If $X$ is $D$-affine, then $H^i(X,D_X)=0$ for $i>0$. An affine variety 
is $D$-affine since $D_X$ is a quasicoherent sheaf of algebras. This implies 
the last assertion.
\end{pf}

\begin{example} Let $X=\bbC^n$. Then $\de (D_X)=H^2(X,\bbC)=0$. 
Since $X$ is affine, $\de (D_X)=\de (D_X(X))$, where $D_X(X)$ is the 
Weyl algebra. It is well known that the Hochschild cohomology of the 
Weyl algebra is trivial.
\end{example}

\section{Deformation of differential operators}

\subsection{Induced deformations of differential operators} 
Let $S$ be a commutative ring and 
$C$ be an $S$-algebra with a finite filtration 
$$0=C_{-1}\subset C_0\subset C_1\subset ...\subset C_n=C,$$
such that the associated graded $grC$ is commutative. Then it makes sense 
to define the ring $D_S(C)=D(C)$ of ($S$-linear) differential operators 
on $C$ in the usual way. More generally, given two left $C$-modules 
$M$, $N$ define the space of differential operators of order $\leq m$ 
from $M$ to $N$ as follows.
$$D^m(M,N)=\{ d\in \Hom _S(M,N)\vert [f_m,...,[f_1,[f_0,d]]...]=0 
\ \text{for all $f_0,...f_m\in C$}\}.$$
Then $D(M,N):=\cup_mD^m(M,N)$ and in particular we obtain a filtered 
(by the order of differential operator) ring $D(C)=D(C,C)$. Note that 
$C\subset D(C)$ acting by left multiplication. Sometimes we will be more 
explicit and will write $D({}_CM,{}_CN)$ for $D(M,N)$. If the algebra 
$C$ is commutative then each $k$-subspace $D^m(M,N)\subset D(M,N)$ is 
also a (left and right) $C$-submodule.

\begin{lemma} Denote by $S_n$ the ring $S[t]/(t^{n+1})$. Then canonically 
$$D_{S_n}(C\otimes _SS_n)\simeq D_S(C)\otimes _SS_n.$$
In particular, for a commutative $k$-algebra $A$ we have 
$$D_{k_n}(A\otimes _kk_n)\simeq D_k(A)\otimes _kk_n.$$
\end{lemma}

\begin{pf} Indeed, every $f\in \End _{S_n}(C\otimes _SS_n)=
\Hom _S(C,C\otimes _SS_n)$ can be uniquely decomposed as 
$$f=\bigoplus_{i=0}^{n}f_i\otimes t^i,$$
where $f_i\in \End _S(C,C)$. Now the inclusion 
$f\in D^m_{S_n}(C\otimes _SS_n)$ is equivalent to inclusions 
$f_i\in D^m_S(C)$ for all $i$. Whence the assertion of the lemma.
\end{pf} 

For the rest of this section we will consider only $k[t]$-algebras, and all 
differential operators will be $k[t]$-linear, so we will omit the 
corresponding subscript. 
We denote as before $k_n=k[t]/(t^{n+1})$. 

Let $A$ be a commutative $k$-algebra and $B$ be a $k_n$-algebra 
with an isomorphism $grB\simeq A\otimes _kk_n$, 
i.e. $B$ defines an element in $\de ^n(A)$. 
Consider the inclusion of rings $D(B)\subset \End_{k_n}(B)$. Both these rings 
are filtered the powers of $t$, hence we obtain a natural homomorphism 
(of degree $0$ of graded algebras).
$$\alpha :grD(B)\to gr\End_{k_n}(B).$$
Note that $\alpha $ may not be injective. 
On the other hand we have a natural homomorphism of graded algebras
$$\delta :gr\End _{k_n}(B)\to \End _{k_n}(grB),$$
which is, in fact, an isomorphism. 

We 
denote the composition of the two maps again by $\gamma :grD(B)\to 
\End _{k_n}(grB)$. 

\begin{lemma} i) The homomorphism $\gamma $ maps $grD(B)$ to 
$D(grB)$.

ii) The following are equivalent

a) The map $\gamma :grD(B)\to D(grB)$ is injective

b) The map $\gamma :grD(B)\to D(grB)$ is surjective.
\end{lemma}

\begin{pf} i). Since everything is $k_n$-linear, it suffices to prove that 
$\gamma (D(B)/tD(B))\subset D(B/tB)$. Let $d\in D^m(B)$ and denote 
by $\bar{d}\in D(B)/tD(B)$ its residue. Let $b_0,...b_m\in B$ with the 
corresponding residues $\bar{b}_0,...,\bar{b}_m\in B/tB$. We have 
$$[b_0,...[b_m,d]...]=0,$$
hence
$$[\bar{b}_0,...[\bar{b}_m,\gamma (\bar{d})]...]=0.$$
Thus $\gamma (\bar{d})\in D^m(B/tB)$.

ii). The injectivity of $\gamma :grD(B)\to D(grB)$ is equivalent to the 
injectivity of the natural map $\alpha :grD(B)\to gr\End _{k_n}(B)$. 
Consider the subspace $D(B/tB)\simeq D(B,t^nB)\subset D(B,B)$. The injectivity 
of $\alpha $ is equivalent to the assertion that every $d\in D(B,t^nB)$ is 
equal to $t^nd_1$ for some $d_1\in D(B)$. But this last assertion is 
equivalent to the surgectivity of the map $D(B)/tD(B)\to D(B/tB)$ and hence to 
the surgectivity if $\gamma :grD(B)\to D(grB)$. 
\end{pf}

\begin{defn} Assume that the map $\gamma :grD(B)\to D(grB)$ is an isomorphism. 
Then by the Lemma 7.1 the algebra $D(B)$ defines an element in $\de ^n(D(A))$. 
We call $D(B)$ the {\it induced} (by $B$) deformation of $D(A)$. We also 
say that $B$ {\it induces} a deformation of $D(A)$. 
\end{defn}

\begin{example} It follows from Lemma 7.1 that the trivial deformation of $A$ 
induces a deformation of $D(A)$, which is also trivial.
\end{example}

\begin{remark} It would be interesting to see which deformations of 
$A$ induce deformations of $D(A)$. 
\end{remark}

\subsection{Two lemmas about induced deformations} 
Assume that $A$ and $B$ are as above and $B$ induces a deformation of $D(A)$. 
Denote the residue  map $\tau :D(B)\to D(A)$. Moreover, assume that 
$D(B)$ is a split extension of $D(A)$ with a splitting homomorphism 
(of $k$-algebras) $s:D(A)\to D(B)$. Since $A\subset D(A)$, the map $s$ 
defines, in particular, a structure of a left $A\otimes _kk_n$-module 
on $B$. The next two lemmas will be used in what follows. 

\begin{lemma} 

i) The residue map $\beta :B\to A$ is a homomorphism of 
left $A$-modules.

ii) $B$ is a free $A\otimes _kk_n$-module of rank 1.
\end{lemma}

\begin{pf} i). 
Given $a\in A$, $b\in B$ we need to show that $\beta (s(a)b)=a\beta (b)$. 
This follows from the identity $\tau s(a)=a$ and the commutativity of 
the diagram
$$\begin{array}{ccc}
D(B)\times B & \stackrel{(\tau ,\beta )}{\longrightarrow} & 
D(A)\times A \\
\downarrow & & \downarrow \\
B & \stackrel{\beta }{\longrightarrow} & A,
\end{array}$$
where the vertical arrows are the action morphisms.

ii) The $A$-module map $\beta :B\to A$ has a splitting $\alpha :A\to B$, 
which induces an isomorphism $\alpha \otimes 1:A\otimes _kk_n\to B$ of left 
$A\otimes _kk_n$-modules.
\end{pf}

\begin{lemma} Assume that the $k$-algebra $A$ is finitely generated.
Consider $B$ with the structure of a left 
$A\otimes _kk_n$-module defined above. Then $D({}_BB)=D({}_{A\otimes _kk_n}B)$ 
as subrings of $\End _{k_n}(B)$.
\end{lemma}

\begin{pf} Denote $\tilde{A}=A\otimes _kk_n$. Since $D(B)$ is a deformation 
of $D(A)$ the graded ring $grD(B)$ coincides with the subring 
$D(grB)\subset \End _{k_n}(grB)$. The isomorphism of $\tilde{A}$-modules 
${}_{\tilde{A}}B\simeq \tilde{A}$ defines an isomorphism of rings 
$$D({}_{\tilde{A}}B)\simeq D(\tilde{A})=D(grB).$$
Hence, in patricular, $grD({}_{\tilde{A}}B)$ is a graded submodule of 
$\End _{k_n}(grB)$ and as such coincides with $D(grB)$. We conclude that 
the graded subrings of $\End _{k_n}(grB)$, $grD(B)$ and $grD({}_{\tilde{A}}B)$ 
coincide $(=D(grB))$. 
So it suffices to prove the inclusion 
$D({}_BB)\subset D({}_{\tilde{A}}B)$. 

We will prove by descending induction on $p$ that 
$$D({}_BB,{}_Bt^pB)\subset D({}_{\tilde{A}}B,{}_{\tilde{A}}t^pB).$$
It follows from Lemma 7.6,i) that the $A$- and $B$-module structure on 
$B$ coincide modulo $t$. More precisely, if $b\in B$ and $a=\beta (b)\in A$, 
then $$s(a)-b:t^\bullet B\to t^{\bullet +1}B.$$
This implies that 
$$D({}_BB,{}_Bt^nB)= D({}_{\tilde{A}}B,{}_{\tilde{A}}t^nB).$$

Suppose that we proved the inclusion $D({}_BB,{}_Bt^{p+1}B)
\subset D({}_{\tilde{A}}B,{}_{\tilde{A}}t^{p+1}B).$ Let $a_1,...a_l$ be 
a set of generators of the algebra $A$. Choose $d\in D^m({}_BB,{}_Bt^pB)$. 
Then the operators
$$d_{i_0...i_m}:=[s(a_{i_0}),...,[s(a_{i_m}),d]...],\quad 
i_j\in \{ 1,...,l\}$$
map $B$ to $t^{p+1}B$. Since $s(a_i)\in D({}_BB)$ also 
$$d_{i_0...i_m}\in D({}_BB,{}_Bt^{p+1}B)\subset  
D({}_{\tilde{A}}B,{}_{\tilde{A}}t^{p+1}B).$$
Thus there exists $N$ such that every $d_{i_0...i_m}\in 
 D^N({}_{\tilde{A}}B,{}_{\tilde{A}}t^{p+1}B).$ Since $A$ is commutative this 
implies that for any $c_1,...,c_m\in A$ 
$$[s(c_1),...,[s(c_m),d]...]\in  D^N({}_{\tilde{A}}B,{}_{\tilde{A}}t^{p+1}B).$$
But then $d\in D^{N+m}({}_{\tilde{A}}B,{}_{\tilde{A}}t^pB).$ Hence 
 $D({}_BB,{}_Bt^pB)
\subset D({}_{\tilde{A}}B,{}_{\tilde{A}}t^pB)$, which completes the 
induction step and proves the lemma. 
\end{pf}

\subsection{Sheafification}
Lat $Y$ be a scheme over $k$, $\cB$ -- a sheaf of $k_n$-algebras on $Y$ 
with an isomorphism of sheaves of $k_n$-algebras $gr\cB\simeq \cO_Y
\otimes _kk_n$, i.e. $\cB$ defines an element in $\de ^n(\cO_Y)$. Then 
using the commutator definition as in 9.1 above we define the sheaf $D(\cB)$ 
of $k_n$-linear differential operators on $\cB$. Thus, in particular, 
$D(\cB)$ is a subsheaf of $\cE nd_{k_n}(\cB)$. In this section all the 
differential operators will be $k[t]$-linear, so we omit the corresponding 
subscript. 

As in the ring case we obtain a natural homomorphism of sheaves of 
graded $k_n$-algebras (which, probably, is neither injective, nor surjective 
in general)
$$\tilde{\gamma}:grD(\cB)\to \cE nd_{k_n}(gr\cB).$$
The following two lemmas are the sheaf versions of Lemmas 9.1 and 9.2 which 
will be used later. The proofs are the same.

\begin{lemma} $D(\cO_Y\otimes _kk_n)=D(\cO_Y)\otimes _kk_n \ \ (=D_Y
\otimes _kk_n).$
\end{lemma}

\begin{lemma} The homomorphism $ \tilde{\gamma}$ maps $grD(\cB)$ to 
$D(gr\cB).$
\end{lemma}

\begin{defn} Assume that $\tilde{\gamma}:grD(\cB)\to D(gr\cB)$ is an 
isomorphism. Then by Lemma 9.8 the sheaf $D(\cB)$ defines an element in 
$\de ^n(D_Y)$. We call $D(\cB)$ the {\it induced} (by $\cB$) deformation 
of $D_Y$ and say that $\cB$ {\it induces} this deformation.
\end{defn}

\subsection{Deformations of differential operators on a flag variety}

\begin{thm} Let $G$ be a complex linear simple simply connected algebraic 
group, $B\subset G$ -- a Borel subgroup, $X=G/B$ -- the corresponding flag 
variety. Then any induced deformation of $D_X$ is trivial. 
\end{thm}

\begin{remark} Since $H^1(X,T_X)=0$ (the variety $X$ is rigid) the only 
deformations of $\cO_X$ are ``purely noncommutative'', i.e. they 
correspond to elements of $H^0(X,\wedge ^2T_X)$. In this respect one may 
ask the following question: Suppose $Y$ is a smooth 
projective variety, $\cB$ -- 
a purely noncommutative deformation of $\cO_Y$. Assume that $\cB$ induces a 
deformation $D(\cB)$ of $D_Y$. Is $D(\cB)$ a trivial deformation of $D_Y$?
\end{remark}

\begin{pf} Assume that a sheaf of $k_n$-algebras $\cB$, 
which represents an element 
in $\de ^n(\cO_X)$, induces a deformation (of order $n$) $D(\cB)$ of $D_X$.
Then for any $m>0$ the sheaf $\cB/t^{m+1}\cB$ induces a deformation of order 
$m$, $D(\cB/t^{m+1}\cB)$, of $D_X$. By induction we may assume that 
$D(\cB/t^n\cB)\simeq D_X\otimes _kk_{n-1}$, i.e. $D(\cB)$ represents an 
element in $\de _0^n(D_X)$. (Recall that $\de ^n_0(D_X)\simeq 
\de (D_X)$.) We need to prove that $D(\cB)$ is the trivial element in 
$\de ^n(D_X)$. For simplicity of notation we assume that $n=1$ (the proof 
in the general  case is the same). 

It is well known that $X$ has an open covering $X=\cup_{w\in W}U_w$, where $W$ 
is the Weyl group of $G$ and $U_w\simeq \bbC ^d$, $d=\dim (X)$. Denote the 
covering $\cU =\{U_w\}$. It follows from  Example 8.3 that $D(\cB)_{U_w}$ 
is the 
trivial deformation of $D_{U_w}$ for each $w\in W$. 

The variety $X$ is $D$-affine ([BB]), thus $\de (D_X)
\simeq H^2(X^{an},\bbC )$ (Corollary 8.2). But $H^2(X^{an},\bbC)=
H^{1,1}(X^{an},\bbC )=Pic(X)\otimes _{\bbZ}\bbC$. Let us describe the 
isomorphism $\sigma :Pic(X)\otimes \bbC \to \de (D_X)$ directly. Let 
$\cL$ be a line bundle on $X$. Then $\cL\vert _{U_w}\simeq \cO_{U_w}$ for 
all $w\in W$. Hence $\cL$ is defined by a \v{C}ech 1-cocycle 
$\{ f_{ij}\in \cO ^*_{U_{w_i}\cap U_{w_j}}\}$. Define derivations 
$$\delta _{ij}:D_{U_{w_i}\cap U_{w_j}}\to D_{U_{w_i}\cap U_{w_j}}$$
by the formula
$$\delta _{ij}(d)=[d,log(f_{ij})].$$
Note that though $log(f_{ij})$ is a multivalued analytic function, $[\cdot ,
log(f_{ij})]$ is a well defined derivation of the ring of differential 
operators and it preserves the algebraic operators. So $\delta _{ij}$ is 
well defined. Using these derivations we define the glueing over $U_{w_1}
\cap U_{w_2}$ of the sheaves $D_{U_{w_i}}\otimes \bbC [t]/(t^2)$ and 
$D_{U_{w_j}}\otimes \bbC [t]/(t^2)$. We denote the corresponding 
global sheaf $\sigma (\cL)$. The map $\sigma :Pic(X)\to \de (D_X)$ is a group 
homomorphism which extends to an isomorphism 
$$\sigma :Pic(X)\otimes \bbC \stackrel{\sim}{\rightarrow} \de (D_X).$$

Let us get back to $D(\cB)\in \de (D_X)$. By the above isomorphism,  
$D(\cB)=\sigma (\cL)$ for some $\cL \in Pic(X)\otimes \bbC$. 
We have $D(\cB)_{U_w}=
D_{U_w}\otimes \bbC [t]/(t^2)$, so that $\cB_{U_w}$ has a structure of a 
$D_{U_w}$-module and, in particular, of an $\cO _{U_w}$-module. 
By (a sheaf version of) Lemma 9.6,ii) 
$\cB _{U_w}\simeq \cO _{U_w}\otimes \bbC [t]/(t^2)$ as an 
$\cO _{U_w}$-module. Since the glueing of different $D(\cB)_{U_w}$'s is by 
means of derivations $[\cdot ,log(f_{ij})]$, it follows that the local 
$\cO _{U_w}$-module structure on $\cB$ agree on the intersections 
$U_{w_i}\cap U_{w_j}$. Hence $\cB$ is an $\cO_X$-module, which fits in the 
short exact sequence of $\cO_X$-modules
$$0\to \cO_X \to \cB \to \cO_X \to 0.$$
Since $\Ext ^1_{\cO_X}(\cO_X,\cO_X)=0$, $\cB=\cO_X\otimes \bbC [t]/(t^2)$. 
Thus $D({}_{\cO_X}\cB)=D_X\otimes \bbC [t]/(t^2)$. But by 
(a sheaf version of) Lemma 9.7  
$D({}_{\cO_X}\cB )=D({}_{\cB}\cB)\ (=D(\cB))$, which proves the theorem.
\end{pf}

\end{document}